\documentclass[12pt]{article}
\usepackage{amscd,amsthm}
\usepackage{latexsym}
\usepackage{amsmath}
\usepackage{amssymb}

\begin{document}
\title{On the largest eigenvalue of a sparse random subgraph of the hypercube
}
\author{Alexander Soshnikov\\
University of
California, Davis 
\\Department of Mathematics\\Davis, CA  95616, USA\\
soshniko@math.ucdavis.edu}
\date{}
\maketitle
\begin{abstract}
We consider a sparse random subraph of the $n$-cube where each edge appears
independently with small probability $\ p(n) 
=O( n^{-1+o(1)}). \ $
In the most interesting regime when $  \ p(n) \ $ is not 
exponentially small we prove that
the largest eigenvalue is $ \ \ \Delta(G)^{1/2} \* (1+o(1))= 
\frac{n \* \log 2}{\log(p^{-1})} \* (1+o(1)) \ 
\ $ almost surely,where $ \ \Delta(G) 
 \ $ is the maximum degree of $ \ G.\  $
If $ \ p(n) \ $ is exponentially small but
not proportional to $ \  \ 2^{-n/k}\* n^{-1}, \ \ \ k=1,2,
\ldots, \ $ then with probability going to one 
$ \ \ \lambda_{\max}(G)= \Bigl ( \Delta(G) \Bigr )^{1/2}= \Bigl ( \bigl 
[\frac{n \* 
\log 2}{\log(p^{-1})-\log n}\bigr ]\Bigr )^{1/2}. \ \ $
If  $ \ \ p(n) \ \ $ is proportional to 
$ \  \ 2^{-n/k}\* n^{-1}, \ \ \ k=1,2,
\ldots, \ $ then with probability going to one 
$ \ \ \Bigl ( \Delta(G) \Bigr )^{1/2} \leq \lambda_{\max}(G) <
\Bigl ( \Delta(G) +1 \Bigr )^{1/2}  \ \ $ and $ \ \ |\Delta(G)-\bigl [ \frac{
n\* \log 2}{\log (p^{-1})-\log n } \bigr ] | \leq 1. \ \ $
\end{abstract}

\section{Introduction and Formulation of Results}

Let $ \ Q^n \ $ be a graph of the $n$-cube consisting of $ \ 2^n \ $ vertices
$ \ \ V=\{ x=(x_1,\ldots, x_n); \ \ x_i \in \{0,1 \}, \ i=1,\ldots,n \} \ \ $
and $ \ \ n\*2^{n-1} \ $ edges $ \ \ E=\{ \ \{x,y\}: \ \ \sum |x_i-y_i|=1 \ \}.
\ $  In this paper we study a random subgraph $ \ \ G(Q^n, p(n)), \ \ $
where each edge appears independently with probability $ \ p(n).\ $ Random
subgraphs of the hypercube were  studied by Burtin [5], Erd\"os and 
Spencer [8], Ajtai, Koml\'os and Szemer\'edi  [1] and Bollob\'as [4], among 
others. In particular it was shown that a giant component emerges shortly after
$ \ p=1/n \ $ ([1]) and the graph becomes connects shortly after $ \ p=1/2 \ $
([5],[8],[4]). Recently the model has become of interest in mathematical 
biology ([7], [14], [15]). In
 this paper we are concerned with the behavior of 
the largest eigenvalues of a sparse random graph  $ \ (
 \  p(n)
\leq n^{-1+o(1)})
. \ \ $

 The adjacency matrix of $ G$  is an $ \ 2^n \times 2^n \ $ matrix $A$ whose 
entries are  either one or zero depending on whether the edge $ \ (x,y) \ $ 
is present in $G$ or not. $A$ is a random real symmetric matrix with the 
eigenvalues denoted by $ \lambda_1 \geq \lambda_2 \geq \ldots 
\geq \lambda_{2^n}.
\ $ It follows from the Perron-Frobenius theorem that the largest eigenvalue
is equal to the spectral norm of $A$, i.e. $ \ \ \lambda_{\max}(G)=\lambda_1
=\Vert
A \Vert = \
\max_{j} \* |\lambda_j|. \ \ $
\medskip

\noindent{\bf Remark 1}  
It is easy to see that for a subgraph of the hypercube, or in general, for 
any bipartite graph, $ \ \ \lambda_k (G)= - \lambda_{|V|-k}(G)
, \ \ \ k=1,2,
\ldots, \ \ $ in particular   $ \ \ |\lambda_{\min}(G)|=\lambda_{\max}(G). 
\ \ $

 Our main result is concerned with the asymptotic behavior of the largest 
(smallest) eigenvalues of $A$.

\medskip

\noindent
{\bf Theorem } {\it Let $ \ \ G(Q^n, p(n)) \ \ $ 
is a random subgraph of the $n$-cube and
$ \ \  p(n) \leq n^{-1+o(1)}.
\ \ $ Then the following statements hold.

(i)
If $ \ \  \gamma^{-n} \ll p(n) \leq n^{-1+o(1)} \ \ $ 
for all $ \ \gamma >1, \ $ 
then $ \ \ \lambda_{\max} =\ \ \bigl (\Delta(G)\bigr )^{1/2} \* (1+o(1))=
\Bigl  ( \frac{n \* \log(2)}{\log(p^{-1})} \Bigr  )^{1/2} \* (1+o(1)) \ \ $
almost surely, where  $ \ \Delta(G) \ $ is the maximum degree of $ \ G \ $.
Also
for any $ \ 0 < \alpha < 1 \ $ there exists some positive constant 
depending on $ \ \alpha \ $ such that with  probability
at least $ \ \ 1- 2^{\alpha \* n} \* 
\exp \Bigl (- \bigl ( 2 \* p \* \log (p^{-1})\bigr )^{-2} \Bigr ) \ \ $
there exist at 
least  $ \ \ 2^{[\alpha \* n]}/(2 \* n^2) \ \ $ eigenvalues greater or equal to
$ \ \  \sqrt{ \frac{ (1-\alpha)\* n \* \log 2 \* (1- 1/\log \log (p^{-1}))}
{\log(p^{-1})}-2}.\ \ \ $

(ii) If  $ \ \ 2^{-n}\* n^{-1} \ll p(n) \leq \gamma^{-n} \ \ $ for some $ \ \ 
\gamma \in (1,2], \ \ $ and
$ \ p(n) \ $ is not proportional to $ \  \ 2^{-n/k}\* n^{-1}, \ \ \ k=2,3,
\ldots , \ $ then with probability going to one 
$ \ \ \lambda_{\max}(G)= \Bigl ( \Delta(G) \Bigr )^{1/2}= \Bigl ( \bigl 
[\frac{n \* 
\log 2}{\log(p^{-1})-\log n}\bigr ]\Bigr )^{1/2}. \ \ $

(iii)If
$ \ p(n) \ $ is  proportional to $ \  \ 2^{-n/k}\* n^{-1}, \ \ $ with some
$ \ \ k=2,3,
\ldots , \ $ then with probability going to one
$ \ \ \Bigl ( \Delta(G) \Bigr )^{1/2} \leq \lambda_{\max}(G) <
\Bigl ( \Delta(G) +1 \Bigr )^{1/2},  \ \ $ and $ \ \ |\Delta(G)-\bigl [ \frac{
n\* \log 2}{\log (p^{-1})-\log n } \bigr ]| \leq 1. \ \ $

(iv) If $ \ \  p(n)\* 2^n\*n \to \nu \ \ $ as $ \ \ n \to \infty, \ \ $
then $ \ \ \lambda_{max}(G) \ \ $ converges in distribution to 
$ \ \ Be(e^{-\nu}).
\ \ $

(v) If $ \ \ p(n)\ll 2^{-n}\* n^{-1}, \ \ $ then with probability going to 
one $ \ \ G \ \ $ is empty and $ \ \ \lambda_{max}(G)=0 . \ \ $
}

\medskip

\noindent{\bf Remark 2}  
$ \  \sqrt{\Delta(G)} \ $ is an obvious lower bound for $ \ \lambda_{\max}(G)
\ $ since $ \ \Vert A \* f \Vert^2
 \geq \Delta(G), \ $ where 
$ \ f \ $ is a delta-function with the support at the vertex $ \ x \ $
of the maximum degree $ \ \deg(x)=\Delta(G). \ $

\medskip

\noindent{\bf Remark 3}  
The result  of the theorem  is similar to a recent result of 
Krivelevich and Sudakov [12] on the largest eigenvalue of
a random subgraph $ \ \ G(n,p) \ \ $ of a complete graph,
who proved that 
$$ \lambda_{\max}(G(n,p))= \max ( \Delta(G(n,p))^{1/2}, \ n\*p)
\*(1+o(1)).$$
To some extent our approach has been influenced by [12].
\medskip

 The results claimed to take place almost surely hold with probability 
one on the product of probability spaces corresponding to $ \ G(Q^n, p(n)),
n=1,2,\ldots. \ \ $ We  use the standard notations
$ \ a_n =\Theta(b_n), \ \ a_n=O(b_n) \ \  $ and $ \ a_n=\Omega(b_n) \ \ $
for $ \ \ a_n >0, b_n >0 \ \ $ as $ \ \ n \to \infty \ $ if there exist 
constants $ \ C_1 \ $ and $ \ C_2 \ $ such that $ \ \ C_1 \* b_n < a_n < C_2 \*
\ b_n, \ \ a_n < C_2 \* b_n, \ \ $ or $ \ \ a_n > C_1 \* b_n \ \ $ 
correspondingly.
The equivalent notations $ \ \ a_n=o(b_n) \ \ $ and $ \ \ a_n \ll b_n \ \ $ 
mean that $ \ \ a_n/b_n \to 0 \ \ $ as $ \ \ n \to \infty .\ $

  The rest of the paper is organized as follows.  Section 2 is devoted to the 
proof of Theorem. Several results of auxiliary nature are collected in Section
3.  The concluding remarks are given in Section 4.

  It is a pleasure to thank Sergey Gavrilets and Janko Gravner for
bringing this problem to my attention and for  useful 
discussions. Research was supported in part by the NSF grant DMS-0103948  .

\section{Proof of Theorem}

We start with the case $ \ \ n^{-\Theta(1)} \leq p(n) \leq n^{-1+o(1)}, \ \ $
i.e. $ \ \ 1 \leq
\liminf \frac{
\log(p^{-1})}{\log n}
 \leq \limsup \frac{\log(p^{-1})}{\log n}
 < +\infty. \ \ $
Let us denote $ \ \max (  e^5\*np, \ \exp(\log n/\log \log n)) \ $ by
$ \ \ r_n. \ \ $
We decompose the set of all vertices $V$ into the disjoint union
$ \ \ V=V_1 \sqcup V_2 \sqcup V_3, \ \ $ where
$ \ \ \ V_1=\{ x\in V : \ \ d(x) \leq r_n \},\ \ $
 $ \ \ \ V_2=\{ x \in V : \ r_n < 
d(x)
\leq  \frac{ n }{\log (p^{-1})}\* r_n^{-2}
 , \ \ \ $
and $ \ \ \ V_3 =\{ x\in V : \ \ d(x) > \frac{ n }{\log (p^{-1})} \* r_n^{-2}
\}. \ \ \ $

 We recall that $ \ d(x) \ $ denotes the degree of the vertex $ x$. Let us 
denote the induced graphs by $ \ \ G_i=G[V_i], \ \ i=1,2,3. \ $ We also denote
by $ \ G_4, \ G_5 \ $ and $ \ G_6 \ $ the bipartite subgraphs consisting of all
edges of $ \ G \ $ between $ \ V_1 \ $ and $ \ V_2, \ \ \ V_1 \ $ and $ \ V_3,
\ $ and $ \ \ V_2 \ \ $ and $ \ \ V_3 \ \ $ correspondingly.

\medskip

\noindent{\bf Lemma 1} {\it $ \ \ 
\lambda_{max}(G) \leq \sum_{i=1}^{6} \ \lambda_{
max}(G_i).$}

\medskip

\noindent{\bf Lemma 2}{\it $$
\lambda_{max}(G_1
) \leq \max \Bigl (\exp \bigl
(\frac {
\log n }{
\log \log n
  } \bigr
), e^5 \*
np \Bigr ), $$
$$
\lambda_{max}(G_4) \leq
 (
\frac{ n}{\log n })
^{1/2} \ \Bigl (
 \max (
\exp \bigl
( \frac{
\log n }{
\log \log n  }\bigr
), \ \ 
e^5 \*
np)  \Bigr )^{-
1/2}.$$
}

\medskip

The proofs of Lemmas 1 and 2 are rather standard (see e.g. [13])
and will be omitted.

\medskip

\noindent{\bf Lemma 3} {\it 
\begin{equation}
\begin{split}
&{ \bf E} \bigl ( \# ( x \in V_2 \sqcup V_3 :
\ \ \ \sum_{y \in V_2 \sqcup V_3 \setminus \{x\}} \ \ (A^2)(x,y) >
\frac{n}{\log(p^{-1})} \frac{ \log(\log(n))}{\log(n)} \bigr
 ) \\ &= 
O \Bigl
( \exp \Bigl
( - \frac{1}{4} \*  n \* \exp \bigl ( \log n  / \log \log n
\bigr )\Bigr
)\Bigr ). 
\end{split}
\end{equation}
}

\medskip

{\bf Proof}

Let us denote $ \ \ \frac{n}{\log(p^{-1})} \frac{ \log(\log(n))}{\log(n)}  \ \ 
$ by $ \ D_n. \ $ We estimate the mathematical expectation in (1) from 
above by
$$ 2^n \* \sum_{m=r_n+1}^n \ {n \choose m} \* p^m  \sum_{s_1 + \ldots + s_m >
D_n} \ \prod_{j=1}^m  \ {n \choose s_j} \* p^{s_j}  \Bigl \{ 
{n-2 \choose r_n-2} \ p^{r_n-2} \Bigr \}^{\frac{1}{2} (s_1+\ldots +s_m)} $$
Indeed, we can choose the vertex $ \ x \ $ in $ \ 2^n \ $ ways. The probability
that the degree of $ \ x \ $ is $ \ m, \ \ r_n<m\leq n, \ $ is $
{n \choose m} \*p^m\*(1-p)^{n-m} \leq {n \choose m} \*p^m. \ \ $
We shall call the vertices whose distance in $ \ G \ $ from $ \ x \ $ is one 
by the vertices of the first generation.  Similarly, we shall call the 
vertices 
whose distance from $ \ x \ $ is two by the vertices of the second generation, 
etc. If $ \  \deg(x)=m, \ $ then there are exactly $ \ m \ $ vertices of the 
first
generation. We denote by $ \ s_j, \ \ j=1,\ldots, m, \ \ $ the number of 
edges that  connect the $j$-th vertex of the first generation to $ \ \ V_2 
\sqcup V_3 \setminus \{x\}. \ \ $ Since $ \ \ \sum_{j=1}^m \ s_j =
\sum_{y \in V_2 \sqcup V_3 \setminus \{x\}} \ (A^2)(x,y), \ \ $ we conclude 
that
$ \  s_1 +\ldots +s_m >D_n. \ \ $ The probability 
$$ \ \ \Pr \bigl ( \bigcap_{j=1}^m \* \{ j {\rm th } \ {\rm vertex} \ 
{\rm  has }
 \ \ s_j \ \ 
{\rm edges} \ {\rm connecting } \ {\rm it} \ {\rm  to } 
\ V_2\sqcup V_3\setminus \{x \} \} \bigr )\ \ $$
is estimated from above by
$ \ \ \prod_1^m \* {n \choose s_j} \* p^{s_j}, \ \ $ since the events are 
independent. The number of the vertices of the second generation
in $ \ \ V_2 \sqcup V_3 \setminus \{x \} \ \ $ is at least
$ \ \ \frac{1}{2} ( s_1+\ldots +s_m) \ \ $ (each vertex of the second 
generation is connected to at most two vertices of the first generation). 
Finally, each of the vertices of the second generation is connected to at least
$ \ \ (r_n-2) \ \ $ vertices of the third generation. Because the edges are 
independent, the last factor in the bound of the mathematical expectation is
$$ \Bigl \{ {n-2 \choose r_n-2}\* p^{r_n-2} \Bigr \}^{\frac{1}{2}\* (s_1+\ldots
s_m)}. $$
Let us denote the sum $ \ \ s_1+ \ldots + s_m \ \ $ by $ \ B \ $ and the term
$$ 2^n \*  {n \choose m} \* p^m   \ \prod_{j=1}^m  \ {n \choose s_j} \* 
p^{s_j}  \Bigl \{ 
{n-2 \choose r_n-2} \ p^{r_n-2} \Bigr \}^{\frac{1}{2} (s_1+\ldots +s_m)} $$
by  $ \ \ T(s_1,\ldots, s_m). \ \ $ We estimate $ \ \ \log T(s_1,\ldots s_m)
\ \ $ from above as 
\begin{equation*}
\begin{split}
& \log T \leq n\*\log 2 - m\* \log m + m \*\log(np) + m + \sum_{j=1}^{m}
( - s_j\* \log s_j + s_j \* \log(np) + s_j) +\\
&\frac{1}{2} \* ( -r_n \* \log r_n
+ r_n \* \log(np) + r_n)(s_1+\ldots +s_m) \leq
n\*\log 2 - \frac{1}{3} r_n \* (s_1+\ldots +s_m).
\end{split}
\end{equation*}

Then the mathematical expectation in (1)
is bounded from above by
\begin{equation*}
\begin{split}
&\sum_{m=r_n+1}^{n} \sum_{B \geq D_n} {B+m-1 \choose m-1}\*\exp( n\*\log 2 - 
\frac{1}{3} \* r_n \* B) = \\ &O \Bigl ( \exp \bigl (-\frac{1}{4} \* n \*
\exp(\log n/\log \log n)\bigr) \Bigr).
\end{split}
\end{equation*}
Lemma 3 is proven.

\medskip

\noindent{\bf Lemma 4} {\it 
\begin{equation}
\begin{split}
& {\bf E} \Bigl( \# \bigl ( x \in V_1 : \ \sum_{y=1}^{2^n} (A^2)(x,y) > \frac 
{ n \* \log 2}{\log ( p^{-1})}\* \bigl (1 + \frac{4}{\log \log n} +\frac{2
\log r_n}
{\log n} \bigr )\bigr ) \Bigr )=\\
& O \Bigl ( \exp \bigl( - \frac {n}{ \log\log n}\bigr )
\Bigr ).
\end{split}
\end{equation}
}

\medskip

{\bf Proof}

Let us denote $ \ \ 
\frac { n \* \log 2}{\log ( p^{-1})}\* (1 + \frac{
4}{\log\log n} + \frac{2 \*\log r_n
}{\log n}
) \ \ $ by $ \ \ L_n.
\ \ $ Let $ \ \ 1 \leq  m=\deg(x) \leq r_n \ \ $ be the number of the vertices 
of the first
generation. Then there are at least $ \ L_n-m \ $ vertices of the second 
generation and similarly to the proof of Lemma 3 we can estimate the l.h.s. of $(2)$
from above by
\begin{equation*}
\begin{split}
& 2^n \sum_{m=1}^{r_n} {n \choose m} \* p^m \* {(n-1)\*m \choose L_n-m} \* 
p^{L_n-m} \leq \\
&r_n \* \max_{1\leq m \leq r_m} \Bigl ( \exp \bigl ( n\*\log 2 - m\*
\log m
 + m
+ m\*\log(np) + (L_n -m)\* \log m  - \\&
(L_n
 -m)\* \log (L_n -m) + 
(L_n-m) \* \log(np)
+ L_n-m \bigr ) \Bigr ) \leq \\
& n \* \exp \Bigl ( n\* \log 2  + r_n \* \log n - L_n \*  \bigl
(\log L_n - 
\log(np) - 
\log r_n - \frac{r_n}{L_n} \* \log L_n  \bigr
) \Bigr )
 \leq \\
& \exp \Bigl ( n \* \log 2 + (r_n +1)\* \log n  - \frac{n \* \log 2}
{\log (p^{-1})}\* \bigl ( 1 +
  \frac{4}{\log
\log n} +\frac{2\* \log r_n}{\log n} \bigr
)\times
\\&
   \bigl
( \log (p^{-1}) + O(1) - 
\log\log ( p^{-1}) - \log r_n \bigr )\Bigr ) \leq \\
& O\Bigl ( - \exp \bigl (\frac {n}{\log \log n} \bigr ) \Bigr ).
\end{split}
\end{equation*}
Lemma 4 is proven.

\medskip

\noindent{\bf Lemma 5} {\it 
\begin{equation}
\begin{split}
&
{\bf E} \Bigl ( \# \bigl ( x \in V_3 : \sum_{y \in V_2 \sqcup V_3} (A(x,y))
^2 >
\frac{n \* \log 2}{\log (p^{-1}) \log n }\bigr )  \Bigr )=\\
& O \Bigl
( \exp \bigl
( -n \* \exp 
( \frac{\log n}{2 \* \log \log n} )\bigr ) \Bigr ).
\end{split}
\end{equation}
}

\medskip

{\bf Proof}

 We estimate the l.h.s. of $ (3)$ form above by 
\begin{equation*}
\begin{split}
&
2^n \sum_{m >\frac{n\* \log 2}{\log (p^{-1}) \* \log n}} {n \choose m} \* p^m 
\Bigl \{ {n \choose r_n}\* p^{r_n} \Bigr \}^m \leq \\
& O \bigl (
(\exp(n\* \log 2 + \log n  - r_n \* \frac{n\* \log 2}{\log (p^{-1}) 
\* \log n})\bigr ) \leq \\
& O \Bigl
( \exp \bigl
( -n \* \exp \bigl
( \frac{\log n}{2 \* \log \log n} \bigr )\bigr ) \Bigr ).
\end{split}
\end{equation*}
Lemma 5 is proven.

\medskip

\noindent{\bf Lemma 6} {\it 
\begin{equation}
\begin{split}
& {\bf E}
 \Bigl ( \# \bigl ( x \in V_3 :  \deg_{G_3}(x)=
\sum_{y \in V_3} (A(x,y))^2 > \bigl (
\frac{n \* \log 2}{\log (p^{-1}) \log n }\bigr )^{1/2}  \bigr ) \Bigr ) =\\
& O \Bigl ( \exp \bigl ( -n^{4/3} ) \Bigr ).
\end{split}
\end{equation}
}
\medskip

{\bf Proof}

Let us denote $ \ \ \frac { n \* \log 2}{ \log (p^{-1})} \* r_n^{-2} \ \ $
by $ \ M_n \ $ and $ \ \ 
\bigl (
\frac{n \* \log 2}{\log (p^{-1}) \log n }\bigr )^{1/2} \ \ $ by $ \ N_n.\ $
We estimate the l.h.s. of (4) from above by
\begin{equation*}
\begin{split}
&
2^n \sum_{m=N_n+1}^n {n \choose m} \* p^m \bigl \{ {n \choose M_n} p^{M_n} 
\bigr \}^m \leq \\
& n \* \exp \Bigl ( n \* \log 2 - \bigl (
M_n \* \log M_n - M_n - M_n \* \log (np) \bigr )\* N_n \Bigr ) = \\
& \exp \Bigl (
   -n^{3/2} \* r_n^{-3} \Bigr )= O \bigl
 ( \exp (-n^{4/3}) \bigr ).
\end{split}
\end{equation*}
Lemma 6 is proven.
\medskip

\noindent{\bf Lemma 7} {\it Let $ \ G \ $ be a random subgraph of the 
$ \ n$-cube, $ \ \ G=G(Q^n,p(n)), \ \ $ where $ \ \ p(n) \leq n^{-1+o(1)}
. \ \  $ Let us define $ \ \ \kappa(n):= \max \{ k:  2^n\* {n \choose k} \*
 p^k \* (1-p)^{n-k} \geq 1 \}. \ \ $ Then the following statements hold.

(i) If $ \ p(n) \ $ is not exponentially small in $ n, $ then
\begin{equation}
\Pr \bigl ( \Delta(G) < \kappa(n) -j \bigr ) \leq
\exp \Bigl ( - \Bigl (\frac{\log 2}{ p \* \log( p^{-1})} 
 \Bigr )^{j} \* (1+o(1)) \Bigr ),
\end{equation}
\begin{equation}
\Pr \bigl ( \Delta(G) > \kappa(n) +j \bigr ) \leq
\Bigl ( \frac{ p \* \log( p^{-1})}{\log 2} \Bigr )^{j} \* (1+o(1)),
\end{equation}
for $ \ j =1,2,\ldots. \ $ 
In particular, with probability one,  there exists sufficiently large (random)
$ \ n_* \ $ such that for $ \ n > n_* \ $ we have
$ \ \ | \Delta(G)- \kappa(n) |\leq 2.$

(ii) If $ \ p(n)= \Theta( 2^{-n/k} \* n^{-1}), \ $ then 
$ \ \  2^n \* { n \choose k } \*p^k =\Theta(1), \ 
 \ \ \kappa(n)= k-1 \  $ or $ \  k,  $ and
\begin{equation}
\begin{split}
& \Pr \bigl ( \Delta(G) =k-1 \bigr )= \exp \bigl ( - 2^n \* { n \choose k } \*
p^k \bigr )\* (1+o(1)),\\
& \Pr \bigl (\Delta(G)=k-1 \bigr )+\Pr \bigl ( \Delta(G) =k \bigr )=
1-O (2^{-n/k}).
\end{split}
\end{equation}
(iii) If $ \ \ p(n) \ \ $ is exponentially small, but not proportional
to $ \ \ 2^{-n/k}\* n^{-1}, \ \ $ then
$ \ \ \kappa(n) = \bigl [ \frac{n \* \log 2}{\log (p^{-1}) - \log n} \bigr ],
\ \  \ \  {\bf E} \* X_{\kappa(n)+1} \ll 1 \ll  
{\bf E} \* X_{\kappa(n)},  \ \ $ 
and
\begin{equation}
\Pr \bigl ( \Delta(G) > \kappa(n) \bigr )= O ( {\bf E} \* X_{\kappa(n) +1)}),
\end{equation}
\begin{equation}
\Pr \bigl ( \Delta(G) < \kappa(n) \bigr ) \leq 
\exp \bigl (- \Theta ( {\bf E} \* X_{\kappa(n)})\bigr ).
\end{equation}
}

\medskip

{\bf Proof}

 Let us denote the number of vertices of $ \ G(Q^n,p) \ $ with degrees larger 
than $ \ k-1 \ $ by $ \ X_k. \ $ Then $ \ X_k= \sum_{i=1}^{2^n} I_i, \ $ where
we denoted by the first $ \ 2^n \ $ positive integers the vertices of $ \ Q^n 
\ $  and by $ \ I_i \ $  the indicator of the event that 
$ \ \deg(i) \geq k. \ $  By its
definition  $ \ X_k \ $ is monotone (non-increasing) with respect to $ \ k. \ $
One can easily calculate the mathematical expectation
$ \  {\bf E} \* X_k = 2^n \sum_{l\geq k } {n \choose l} p^l (1-p)^{n-l}.\  $
Estimating $ \ {\bf E} \*
 X_k \ $ from above as
$$ {\bf E} \* X_k \leq 2^n \* {n \choose k} \* p^k  \leq 
\exp( n \* \log 2 - k \* 
\log k + k \* \log (n\*p) +O(k) )   $$ we obtain that for
$ \ \ k \geq \frac{ n \* \log 2 \* ( 1+ 1/\log \log (p^{-1}))}
{\log ( p^{-1})} $

\begin{equation}
{\bf E} \* X_{k}=
O\Bigl (\exp \bigl ( -\frac{n}{2 \* \log \log (p^{-1})}\bigr ) 
\Bigr ).
\end{equation}

 On the other hand if  $ \ \ k \leq 
\frac{ n \* \log 2 \* ( 1- 1/\log \log (p^{-1}))}
{\log ( p^{-1})} \ \ $ then 

\begin{equation}
{\bf  E} \* X_{k} = \Omega \Bigl ( \exp \bigl ( \frac{n}
{2 \* \log \log (p^{-1})}\bigr ) \Bigr ). 
\end{equation}


It is  clear that $ \ \kappa(n) \ \ $ must satisfy the inequalities

\begin{equation}
\frac{ n \* \log 2 \* ( 1- 1/\log \log (p^{-1}))}{\log (p^{-1})}-1
\leq \kappa(n) \leq 
\frac{ n \* \log 2 \* ( 1+ 1/\log \log (p^{-1}))}{\log (p^{-1})} +1.
\end{equation}
We claim that for such $ \ k \ $ the probability
$ \ \ \Pr ( \Delta(G) <k) =\Pr (X_k=0)  \ \ $ is equal, up to a  small 
error term, to $ \ \ \exp(-{\bf E} \* X_k). \ \ $ More precisely the following
inequalities take place
\begin{equation}
\begin{split}
& \exp \Bigl (- \frac{{\bf E} \* X_k}{ 1- 2^{-n} \* {\bf E} \*X_k} \Bigr ) 
\leq 
\Pr ( X_k =0)
\leq \\
& \exp \Bigl (- {\bf E} \* X_k \* \bigl ( 1-
{\bf E} \* X_k   \* \frac{p^{-1}}{2^n} \* (\frac{ k^2}{n^2} +p) \* 
\exp( n \* 2^{-n+1} \* {\bf E} \* X_k )\bigr ) \Bigr ). 
\end{split}
\end{equation}
The l.h.s of (13) follows from the FKG inequality ([3],[9]). Since the events
$ \ \ \{ \deg(i)<k \}_{i=1}^{2^n} \ \ $ are monotone with respect to the edge
indicators we have
\begin{equation}
\begin{split}
&\Pr (X_k=0)= \Pr ( \bigcap_{i=1}^{2^n} \{ \deg(i) < k \} ) \geq
\prod_{i=1}^{2^n} (1-{\bf E} I_i) \leq \\
&\prod_{i=1}^{2^n} \exp \Bigl (- \frac{{\bf E} I_i}
{1-{\bf E}\* I_i}\Bigr )  =\exp \Bigl (- \frac{\sum {\bf E} I_i}
{1-\max_{1\leq i \leq 2^n}
{\bf  E} I_i}\Bigr ) =\\
& \exp \Bigl ( -\frac{ {\bf E} \* X_k}{1- {\bf E} \* I_1} \Bigr ).
\end{split}
\end{equation}
The l.h.s. of (13) now follows from $ \ \ {\bf E} I_1  = 2^{-n} \* {\bf E}
 \* X_k. \ \ $
To prove the r.h.s. of (13) we apply the Suen's type inequality ( see e.g.
[9], Theorem 2.22,
part (i)) that states that
\begin{equation}
\Pr(X_k=0)\leq \exp \Bigl ( -{\bf E} \* X_k + \epsilon \* e^{2 \delta} \Bigr ),
\end{equation}
where 
$ \ \ \epsilon= \frac{1}{2} \sum \sum_{i \sim j} {\bf E} (I_i\* I_j), \ \ $
and 
$ \ \ \delta= \max_{i} \sum_{k \sim i} { \bf E} \* I_k. \ \ $
Here we use the notation $ \ i \sim j \ $ if $ \ i \neq j \ $ and $ \ I_i \ $
and $ \ I_j \ $ are dependent random variables.
It is easy to see that
in our case $ \ \ \delta= n\* 2^{-n} \* {\bf E }\* X_k,\ \ $
and
\begin{equation}
\begin{split}
& \epsilon \leq 2^{n-1} \* n \* \Bigl ( p \* \bigl [
{n-1 \choose k-1} p^{k-1}\bigr ]^2
+(1-p) \* \bigl [ {n-1 \choose k} \* p^k \bigr ]^2 \Bigr ) \leq\\
& n\* 2^{n-1} \* 2^{-2\*n} \* (\frac{k^2}{n^2 \* p} +1)  \* ({\bf E} \* X_k)^2
\leq ({\bf E} \* X_k)^2 \* \frac{1}{2^n \* p }\*(\frac{ k^2}{n^2}+p).
\end{split}
\end{equation}

Let us now consider the case (i) in more detail.
Taking into account that $ \ \ {\bf E}  \* X_{k+1} = 
\frac{p \* \log (p^{-1})}{\log 2}
\* {\bf E} \* X_k \* (1+o(1)) \ \ $ for $ \ \ k =\kappa(n) \* (1+o(1)) \ \ $ 
we obtain from the definition of $ \ \ \kappa(n) \ \ $ that for any fixed
$ \ j= 0, 1,2,\ldots, \ $ 
\begin{equation}
 \bigl (\frac{\log 2}{p \* \log (p^{-1})}\bigr )^j \* (1+o(1)) 
\leq {\bf E}  \* X_{\kappa(n)-j} \leq
 \bigl (\frac{\log 2}{p \* \log (p^{-1})}\bigr )^{j+1} \* (1+o(1)),
\end{equation}
\begin{equation}
{\bf E}  \* X_{\kappa(n)+j} \leq
\bigl (\frac{p \* \log (p^{-1})}{\log 2}\bigr )^{j-1} \* (1+o(1)),
\end{equation}
Applying (17) and the r.h.s. of (13)  we infer
\begin{equation}
\begin{split}
&\Pr ( \Delta(G) < \kappa(n)-j) = \Pr (X_{\kappa(n)-j} =0) \leq \\
& \exp(-{\bf E} \* X_{\kappa(n)-j}\* (1+o(1))) \leq  \exp 
\Bigl ( - \bigl ( \frac{\log 2}
{p \* \log(p^{-1}) })^j \* (1+o(1)\bigr ) \Bigr ) \ \ \ j=1,2,..
\end{split}
\end{equation}
In a similar manner
\begin{equation}
\begin{split}
&\Pr ( \Delta(G) > \kappa(n)+j ) = 1- \Pr (X_{\kappa(n)+j+1} =0) \leq \\
& 1- \exp(-{\bf E} \* X_{\kappa(n)+j+1} \* (1+ o(1))) 
\leq  1 - \exp \Bigl ( - \bigl (\frac{ p \* \log(p^{-1})}{\log 2}\bigr )^j 
\* (1+o(1)) \Bigr )\\
& \leq \Bigl (\frac{p\*\log(p^{-1})}{\log 2} \Bigr )^j \* (1+o(1)),
 \ \ \ j=1,2,...
\end{split}
\end{equation}

Let us now consider the case (ii).
Since $ \ \ p(n)= \Theta( 2^{-n/k} \* n^{-1}) \ \ $ we have
$ \ \  {\bf E}  \* X_k= 2^n \* { n \choose k } \*p^k \* (1-p)^{n-k} =
\Theta(1), \ 
 \ $ and $ \ \kappa(n)= k-1 \  $ or $ \  k,  $ depending on whether
$ \ \ {\bf E}  \* X_k < 1 \ $ or 
$ \ \ {\bf E}  \* X_k \geq 1. \ $ 
It follows from $ \ (13) \ $ that 
\begin{equation}
\begin{split}
& \Pr \bigl ( \Delta(G) < k \bigr )= 
\exp \bigl ( - {\bf E} \* X_k \bigr )\* (1+o(1)) = \\
&\exp \bigl ( - 2^n \* { n \choose k } \*
p^k \bigr )\* (1+o(1)).
\end{split}
\end{equation}
Applying the r.h.s. of $ \ (13) \ $ and 
\begin{equation}
{\bf  E} \* X_{k-1}   \* \frac{p^{-1}}{2^n} \* (\frac{ k^2}{n^2} +p) \* 
\exp( n \* 2^{-n+1} \* {\bf E} \* X_k ) =o(1), \ \ \ \ k\geq 2,
\end{equation}
we obtain
\begin{equation}
\begin{split}
&\Pr \bigl ( \Delta(G) < k-1 \bigr )= 
\exp \bigl ( - {\bf E} \* X_{k-1} \* (1+o(1)) \bigr ) = \\
&  \exp \Bigl ( - \Theta \bigl (\frac
{1}{p \* \log (p^{-1})} \bigr ) \Bigr ) =
\exp \Bigl ( - \Theta \bigl (
2^{-n/k} \bigr ) \Bigr ).
\end{split}
\end{equation}
To estimate
$ \ \ \Pr ( \Delta(G) >k ) \ \ $ we observe that $ \ \ {\bf E} \* X_{k+1} = 
\Theta ( 2^{-n/k}) \ \ $ which implies
\begin{equation}
\Pr ( \Delta(G) >k )= 1- \Pr (X_{k+1}=0) = \Theta (2^{-n/k}).
\end{equation}

 If $ \ \ p(n) \ \ $ is exponentially small but not proportional
to $ \ \ 2^{-n/k}\* n^{-1},  \ \ k=1,2,\ldots, $ then
$ \ \ \kappa(n) = \bigl [ \frac{n \* \log 2}{\log (p^{-1}) - \log n} \bigr ]
\ \ $ and $  \ \  {\bf E} \* X_{\kappa(n)+1} \ll 1 \ll  
{\bf E} \* X_{\kappa(n)}.  \ \ $ 
In a similar way to (i), (ii) one has
\begin{equation}
\Pr \bigl ( \Delta(G) > \kappa(n) \bigr )= O ( {\bf E} \* X_{\kappa(n) +1}),
\end{equation}
To estimate
$ \ \ \Pr \bigl ( \Delta(G) < \kappa(n) \bigr )  \ \ $ we consider
first the case \\
$ \ \ p(n) >  2^{-\frac{n}{2}+2}/(n \* \sqrt{\log n}). \ \ $ Then
\begin{equation}
{\bf  E} \* X_{\kappa(n)}   \* \frac{p^{-1}}{2^n} \* 
(\frac{ \kappa(n)^2}{n^2} +p) 
\* 
\exp( n \* 2^{-n+1} \* {\bf E} \* X_{\kappa(n)} ) =o(1), 
\end{equation}
and
\begin{equation}
\Pr \bigl ( \Delta(G) < \kappa(n) \bigr ) \leq 
\exp \bigl (- \Theta ( {\bf E} \* X_{\kappa(n)})\bigr ).
\end{equation}
In the case
$ \ \ 2^{-n} \* n^{-1} \ll p(n) 
\leq   2^{-\frac{n}{2}+2}/(n \* \sqrt{\log n}) \ \ $ one has
$ \ \ \kappa(n) =1 \ \ $ and
$ \ \ \Delta(G) < \kappa(n) \ \ $ iff the graph is empty. This probability
is equal to $ \ \exp \bigl ( - \Theta ( n \* 2^n \* p ) ) \ $
since it is the probability
that $ \ \ n \* 2^{n-1} \ \ $ independent Bernoulli random variables
$ \ \ Be ( p) \ \ $ all equal zero.

Lemma 7 is proven.

Combining the results of Lemmas 1-7 we are now ready to prove part (i) of the 
theorem.
Indeed, applying 
Borel-Contelli Lemma we obtain that with probability one  there exists 
sufficiently large (random) $ \ n_* \ $ such that for all $ \ \ n> n_* \ \ $
the counting numbers from Lemmas 3-6 are all zero. 
Let $ \ \ n> n_*.\ \ $
It follows 
from Lemma 3  that 
\begin{equation}
\begin{split}
& \max_{x \in V_2} \sum_{y \in V_2} (A^2)(x,y) \leq \frac{n}{\log(p^{-1})}\*
\frac{\log \log n}{\log n} + \max_{x \in V_2} (A^2)(x,x) =\\
&\frac{n}{\log(p^{-1})}\*
\frac{\log \log n}{\log n} + \max_{x \in V_2} \deg_{G_2}(x) \leq \\
&\frac{n}{\log(p^{-1})}\*
\frac{\log \log n}{\log n} + \frac{n \* \log 2}{\log(p^{-1})}\* r_n^{-2}=\\
&o\bigl (\frac{n}{\log(p^{-1})}\bigr ).
\end{split}
\end{equation}
Since
$$ \ \  \bigl (
\lambda_{\max}(G_2) \bigr )^2 \leq \max_{x \in V_2} \sum_{y \in V_2} (A^2)(x,y)
\ \ $$we conclude that
\begin{equation}
\lambda_{\max}(G_2)=o \Bigl (
\bigl(\frac{n}{\log(p^{-1})}\bigr )^{1/2} \Bigr )
\end{equation}
almost surely. Similar estimates hold for $ \ \ \lambda_{\max}(G_3) \ \ $
and $ \ \ \lambda_{\max}(G_6). \ \ $
The estimate
\begin{equation}
\lambda_{\max}(G_3)=o \Bigl ( \bigl
(\frac{n}{\log(p^{-1})} \bigr )^{1/2} \Bigr )
\end{equation}
follows from 
$  \ \ 
\lambda_{\max}(G_3) \leq \max_{x \in V_3} \deg_{G_3}(x) \ \ $
and Lemma 6. To prove 
\begin{equation}
\lambda_{\max}(G_6)=o \Bigl ( \bigl
(\frac{n}{\log(p^{-1})} \bigr )^{1/2} \Bigr )
\end{equation}
we employ $(28)$, Lemma 3 and Lemma 5 to see that
\begin{equation}
\begin{split}
&\max_{x \in V_3} \sum_{y \in V_2, z \in V_3} A(x,y)\* A(y,z)=
o(\frac{n}{\log(p^{-1})}), \\
&
\max_{x \in V_2} \sum_{y \in V_3, z \in V_2} A(x,y)\* A(y,z)=
o(\frac{n}{\log(p^{-1})}).
\end{split}
\end{equation}
Finally, we claim that
\begin{equation}
\lambda_{\max}(G_5)=\Bigl ( \bigl
(\frac{n\*\log 2}
{\log(p^{-1})} \bigr )^{1/2} \Bigr) \  \ {\rm (a.s)},
\end{equation}
which
follows from Lemmas
3,4 and 7. Combining Lemmas 1, 2 and $ (29)-
(32) \ $ we prove
$ \ \ \lambda_{\max}(G(Q^n,p))= \bigl ( \Delta(G(Q^n,p))\bigr )^{1/2}\* 
(1+o(1))
\ \ $ almost surely for
$ \ \ n^{-\Theta(1)} \leq p(n) \leq n^{-1+o(1)}. \ \ $
To find the the eigenvalues of $ \ A \ $ close to the 
$ \ \lambda_{\max} \ $ we use Lemmas 8 and 9 from the 
next section. Assuming Lemma 8 we can construct 
$ \ 2^{[\alpha \* n ]} \  \ \delta$-functions $ \ \{ f_i\}_{i=1}^{2^{[\alpha 
\* n ]}}\ \ $  supported at the   vertices$ \ \ \{x_i\}_{i=1}^{2^{[\alpha
\* n]}}\ \ $ of  degrees greater or equal than $ \ \kappa (n -[\alpha \* n]) -2
\geq \frac{(1-\alpha) \* n \* \log 2 \* (1-1/\log \log (p^{-1})) }
{\log (p^{-1})}-2. \ \ $
The constructed vectors
$ \ \{ f_i\}_{i=1}^{2^{[\alpha \* n]}} \ \ $ form an orthonormal family
such that
$ \ \ \bigl ( A^2 \* f_i, f_i \bigr ) >   \ \ \frac{(1-\alpha) \* n \* 
\log 2}{\log (p^{-1})}.  \ \ $ Since each vertex has at most $ \ \ (n^2-1) \ \ $
vertices within distance 2 one can select a sub-family of size at least 
$ \ \ \frac{2^{[\alpha \* n]}}{n^2} \ \ $ such that for the  vectors from the 
sub-family
$  \ \ \bigl ( A^2 \* f_i, f_j \bigr ) =0  \ \ $ if $ \ \ i\neq j . \ \ $
Applying Lemma 9 one obtains that there are at least 
$ \ \ \frac{2^{[\alpha \* n]}}{n^2} \ \ $  eigenvalues of $ \ A^2 \ $ 
greater or equal to $ \ \kappa (n -[\alpha \* n]) -2. \ $
Since the spectrum
of $ \ A \ $ is central symmetric with respect to the origin this implies
that there are at least $ \ \ \frac{2^{[\alpha \* n]}}{2 \*
n^2} \ \ $ 
 eigenvalues of $ \ A   \ $  greater or equal to
$ \ \ \Bigl (\kappa (n -[\alpha \* n]) -2
 \Bigr )^{1/2}
. \ \ $ 

  The case $ \ \ \log n \ll \log (p^{-1}) \ll n \ \ $ is very similar
to the previous one. We again represent $ \ \ V \ \ $ as
$ \ \ V=V_1 \sqcup V_2 \sqcup V_3, \ \ $ where now
$ \ \ \ V_1=\{ x\in V : \ \ d(x) \leq \tau_n \},\ \ $
 $ \ \ \ V_2=\{ x \in V : \ \tau_n < 
d(x)
\leq  \frac{ n }{\log (p^{-1})}\* \tau_n^{-2}
 , \ \ \ $
and $ \ \ \ V_3 =\{ x\in V : \ \ d(x) > \frac{ n }{\log (p^{-1})} \* 
\tau_n^{-2},
 \ \ \ $ with $ \ \ \tau_n= \exp(\frac{\log(\Delta_*)}{\log\log(\Delta_*)}),
\ \ $ and $ \ \ \Delta_* = \frac{n \* \log 2}{\log (p^{-1})}. \ \ $
We claim that the analogues of Lemmas 3-6 hold, namely:

\begin{equation}
\begin{split}
&{\bf E} \bigl ( \# ( x \in V_2 \sqcup V_3 :
\ \ \ \sum_{y \in V_2 \sqcup V_3 \setminus \{x\}} \ \ (A^2)(x,y) >
\frac{n}{\log(p^{-1})\* \tau_n^{1/3}} \bigr
 ) \\ &= 
O \bigl
( \exp \bigl
( - n \* \tau_n^{1/2}\bigr
)\bigr ),
\end{split}
\end{equation}

\begin{equation}
\begin{split}
& {\bf E} \bigl( \# ( x \in V_1 : \ \sum_{y=1}^{2^n} (A^2)(x,y) > \frac 
{ n \* \log 2}{\log ( p^{-1})}\* (1 +1/\log\log(p^{-1}))
) \bigr )=\\
& O \bigl ( \exp \bigl( - \frac {n}{ \log\log (p^{-1})}\bigr )
\bigr ),
\end{split}
\end{equation}

\begin{equation}
{\bf E} \Bigl ( \# \bigl ( x \in V_3 : \sum_{y \in V_2 \sqcup V_3} (A(x,y))
^2 >
\frac{n \* \log 2}{\log (p^{-1}) \tau_n^{1/2} }\bigr )  \Bigr )=
 O \Bigl
( \exp \bigl
( -n \* \tau_n^{1/2}\bigr ) \Bigr ),
\end{equation}
and
\begin{equation}
\begin{split}
& {\bf E}
 \Bigl ( \# \Bigl ( x \in V_3 :  \deg_{G_3}(x)=
\sum_{y \in V_3} (A(x,y))^2 > \Bigl (
\frac{n \* \log 2}{\log (p^{-1}) \tau_n }\Bigr )^{1/2}  \Bigr ) \Bigr ) =\\
& O \Bigl ( \exp \bigl ( -n \* \tau_n ) \Bigr ).
\end{split}
\end{equation}

The proofs of $ \ \ (34)-(37) \ \ $ are very similar to the arguments given in 
Lemmas 3-6 and  left to the reader.

 Let us now consider the case when $ \ p(n) \ $ is exponentially small in $ \ 
n . \ $  We denote by $ \ Y_k \ $ the number of isolated components with
$ \ k \ $ edges, $ \ \ k=1,2,\ldots. \ \ $  It is easy to see that 
\begin{equation}
{\bf E} \* Y_k =\Theta ( {\bf E} \* X_k )= \Theta ( 2^n \* n^k \* p^k ).
\end{equation}
 If $ \ \ p(n) \ \ $ is
not proportional to  $ \  \ 2^{-n/k}\* n^{-1}, \ \ \ k=1,2,3,
\ldots, \ \ $ then it follows from Lemma 7, part (iii) and (38) that with
probability going to one the maximum degree of $ \ G(Q^n,p) \ $ is 
$ \ \kappa(n)=\bigl [ \frac{n \* \log 2}{\log (p^{-1}) - \log n} \bigr ]
\ $ and there are no components with more than $ \ \kappa(n) \ $ edges.
Since the largest eigenvalue of $ \ G \ $ is the maximum of the eigenvalues of 
its connected components and  the largest eigenvalue of a component with $ \ 
k \ $ edges is not greater than $ \ \sqrt{k} \ $ ( and is equal to 
$ \ \sqrt{k} \ $ only if the component is a star on $ \ k+1 \ $ vertices),
we prove that with probability going to one
$ \ \lambda_{\max}(G)=\sqrt{\kappa(n)}=\sqrt{\Delta(G)}. \ \ $

  Finally if  $ \ \ p(n) \ \ $ is
proportional to  $ \  \ 2^{-n/k}\* n^{-1}, \ \ \ k=1,2,3,
\ldots, \ \ $ then with probability going to one 
$ \ \Delta(G)  \in \{ k-1,k\} \ $ and there are no
connected components with more than $ \ k \ $ edges.

Theorem is proven.

\medskip

\section{Auxiliary Results}

In this section we present two auxiliary lemmas.
Our first result claims that there are many vertices with degrees close 
to the maximum degree

\medskip

\noindent{\bf Lemma 8} {\it 
If $ \ p(n) \ $ is not exponentially small,
then  for any fixed $ \ \ 0 < \alpha < 1  \ \ $
with probability at least
$ \ \ 1- 2^{\alpha \* n} \* 
\exp \Bigl (- \bigl ( 2 \* p \* \log (p^{-1})\bigr )^{-2} \Bigr ) \ \ $
there exist at 
least  $ \ \ 2^{[\alpha \* n]}/(2 \* n^2) \ \ $ vertices with degrees 
 greater or equal to $ \ \ 
\frac{ (1-\alpha)\* n \* \log 2 \* 
(1- 1/\log \log (p^{-1}))}{\log(p^{-1})} -2.\ \ \ $
}

\medskip

{\bf Proof}

Consider $ \ \ Q^n \ \ $ as a disjoint union of $ \ \ 2^{[\alpha\*n]} \ \ $
cubes of dimension $ \ \ n- [\alpha \* n], \ \ \ 
Q^n= \sqcup_{i=1}^{2^{[\alpha\*n]}} Q_i. \ \ $ Consider
random subgraphs
$ \ \ G_i= G(Q_i, p(n)) \ \ $ induced by the edges of $ \ G. \ $
According to Lemma 7 the maximum degree of $ \ \ Q_i  \ \ $ is at least
$ \ \ \kappa( n -[ \alpha \* n ] )-2 \geq
  \frac{ (1-\alpha)\* n \* \log 2 \* 
(1- 1/ \log \log (p^{-1}))}{\log(p^{-1})}-2 \ \ \ $ with probability
at least $ \ \ 
\exp \Bigl (- \bigl ( 2 \* p \* \log (p^{-1})\bigr )^{-2} \Bigr ). \ \ $
The intersection of these events has probability at least
$ \ \ 1- 2^{\alpha \* n} \* 
\exp \Bigl (- \bigl ( 2 \* p \* \log (p^{-1})\bigr )^{-2} \Bigr ). \ \ $

Lemma 8 is proven.

We finish this section with elementary lemma from linear algebra.

\medskip

\noindent{\bf Lemma 9} {\it Let $ \ \ A \ \ $ be a Hermitian 
(or real symmetric) matrix and $ \ \ \{ f_i \}, \ i=1,\ldots n, \ \ $ be an 
orthonormal family of vectors such that $ \ \ (A \* f_i, f_i ) \geq \lambda
 \ \ $ for some $ \ \ \lambda, \ \ i=1,2,\ldots,n. \ $ Suppose that
$ \ \ (Af_i,f_j)=0 \ \ $ if $ \ i \neq j. \ $ Then the number of eigenvalues 
of $ \ A \ $  greater or equal to $ \ \lambda \ $ is at least $ \ n. \ $}

\medskip

{\bf Proof}

Let the number of the eigenvalues greater or equal to $ \ \lambda \ $ be less 
than $ \ n. \ $ Then there exists a non-zero linear combination
$ \ \ f= \sum_{i=1}^n x_i \* f_i \ \ $ orthogonal to all eigenvectors
with the eigenvalues greater or equal to $ \ \lambda, \ $ which implies
$ \ \ (Af,f) < \lambda \* (f,f). \ \ $ On the other hand it follows from the 
conditions of the lemma that
\begin{equation*}
(Af,f) = \sum_{i=1}^n |x_i|^2 \* (Af_i,f_i) \leq 
\lambda \* 
\sum_{i=1}^n |x_i|^2 = \lambda \* (f,f).  
\end{equation*}

\medskip

\section{Concluding Remarks}

\medskip

It would be also interesting to study  the regime
$ \ \ n^{-1+o(1)} \leq p(n) \leq 1 \ \ $ and to prove the analogue of the 
Krivelevich- Sudakov theorem there as well.

 There are several other important questions that are beyond the reach of 
presented technique. The most fundamental is perhaps the local statistics
of the eigenvalues, in particular the local statistics near the edge of the 
spectrum. For the results in this direction for other random matrix models
we refer the reader  to [18],[19],[17].  A recent result of Alon, Krivelevich 
and Vu [2] states that the deviation of the first, second, etc largest 
eigenvalue from its mean is at most of order of $ \ O(1). \ $ Unfortunately
our results  give only the leading term of 
the mean.

  Second, and perhaps even more difficult question  is whether the local 
behavior of the eigenvalues is not sensitive to the details of the distribution
of the matrix entries of $ \ A. \ $  We refer the reader to
[16],[6],[17],[10] for the results of that nature for unitary invariant and 
Wigner random matrices.

\def\cmp{{\it Commun. Math. Phys.} }

\end{document}